\documentclass[12pt]{article}

\usepackage{hyperref}
\usepackage{amsmath}
\usepackage{amssymb}
\usepackage{amsfonts}
\usepackage{amsopn}
\usepackage{amsthm}
\usepackage[all]{xy}

\swapnumbers
\theoremstyle{plain}
\newtheorem{thm}{Theorem}[section]
\newtheorem{lem}[thm]{Lemma}
\newtheorem{cor}[thm]{Corollary}
\newtheorem{prop}[thm]{Proposition}

\theoremstyle{definition}
\newtheorem{ntt}[thm]{}

\newtheorem{rem}[thm]{Remark}

\newcommand{\pl}{\mathbb{P}}   
\newcommand{\zz}{\mathbb{Z}}   
\newcommand{\A}{\mathrm{A}} 
\newcommand{\B}{\mathrm{B}} 
\newcommand{\C}{\mathrm{C}} 
\newcommand{\FF}{\mathrm{F}_4} 
\newcommand{\G}{\mathrm{G}_2}  
\newcommand{\id}{\mathrm{id}}       
\newcommand{\pr}{\mathrm{pr}}       
\newcommand{\can}{\mathrm{can}}     
\newcommand{\res}{\mathrm{res}}     
\newcommand{\op}{\mathrm{op}}     

\newcommand{\Var}{\mathcal{V}ar} 
\newcommand{\M}{\mathcal{M}}     
\newcommand{\J}{\mathcal{J}} 
\newcommand{\E}{\mathcal{E}} 
\newcommand{\F}{\mathcal{F}} 
\newcommand{\V}{\mathcal{V}} 
\newcommand{\OO}{\mathcal{O}} 
\newcommand{\Cor}{\mathcal{C}or} 
\newcommand{\zab}{\mathbb{Z}\text{-}\mathcal{A}b} 
\newcommand{\Sb}{\mathcal{H}} 
\newcommand{\De}{\mathit{\Delta}} 
\newcommand{\Seg}{\mathrm{Seg}} 

\newcommand{\xra}[1]{\xrightarrow{#1}} 

\DeclareMathOperator{\Br}{\mathrm{Br}}     
\DeclareMathOperator{\CH}{\mathrm{CH}}      
\DeclareMathOperator{\Aut}{\mathrm{Aut}}    
\DeclareMathOperator{\Hom}{\mathrm{Hom}}    
\DeclareMathOperator{\Mor}{\mathrm{Mor}}    
\DeclareMathOperator{\End}{\mathrm{End}}    
\DeclareMathOperator{\Coker}{\mathrm{Coker}}    
\DeclareMathOperator{\im}{\mathrm{im}}      
\DeclareMathOperator{\SB}{\mathrm{SB}}      
\DeclareMathOperator{\PGL}{\mathrm{PGL}}      
\DeclareMathOperator{\Gr}{\mathrm{Gr}} 
\DeclareMathOperator{\rdim}{\mathrm{rdim}} 
\DeclareMathOperator{\rk}{\mathrm{rk}} 
\DeclareMathOperator{\tr}{\mathrm{tr}} 
\DeclareMathOperator{\ind}{\mathrm{ind}} 

\title{Chow motives of twisted flag varieties\footnote{MSC2000: 57T15, 19E15}}

\author{B.Calm\`es, V.Petrov, N.Semenov, K.Zainoulline}

\begin{document}

\maketitle

\begin{abstract}
Let $G$ be an adjoint simple algebraic group of inner type.
We express the Chow motive (with integral coefficients)
of some anisotropic projective $G$-homogeneous
varieties in terms of motives of simpler $G$-homogeneous varieties,
namely, those that correspond to maximal parabolic subgroups
of $G$.
We decompose the motive of a generalized Severi-Brauer variety $\SB_2(A)$,
where $A$ is a division algebra of degree 5, into a direct
sum of two indecomposable motives. As an application we provide
another counter-example to the uniqueness of a direct sum decomposition in
the category of motives with integral coefficients.\\

Keywords: anisotropic projective homogeneous variety, Chow motive
\end{abstract}

\section{Introduction}

Let $G$ be an adjoint simple algebraic group of inner type
over a field $k$. Let $X$ be a twisted
flag variety, i.e., a projective $G$-homogeneous variety over $k$.
The main purpose of the paper is to express
the Chow motive of $X$ in terms of motives of ``minimal'' flags, i.e.,
those $G$-homogeneous varieties
that correspond to maximal parabolic subgroups of $G$.

Observe that the motive of an \emph{isotropic} $G$-homogeneous variety
can be decomposed in terms of motives
of simpler $G$-homogeneous varieties using the techniques
developed by Chernousov, Gille, Merkurjev \cite{CGM} and
Karpenko \cite{Ka01}. For $G$-varieties, when $G$ is \emph{isotropic},
one obtains a similar decomposition
following the arguments of Brosnan \cite{Br03}.
In the case of $G$-varieties, where $G$ is {\it anisotropic},
no general decomposition methods are known
except several particular cases of quadrics (see for example Rost \cite{Ro98})
Severi-Brauer varieties
(see Karpenko \cite{Ka95}) and
exceptional varieties of type $\FF$
(see \cite{NSZ}).

In the present paper we provide methods that allow to decompose
the motives of some \emph{anisotropic}
twisted flag $G$-varieties,
where the root system of $G$ is of
types $\A_n$,
$\B_n$, $\C_n$, $\G$ and $\FF$,
i.e., has a Dynkin diagram which does not branch.

As an application,
we provide another counter-example to the uniqueness
of a direct sum decomposition in the category of Chow motives
with integral coefficients (see \ref{cex}).
Observe that such a counter-example was already constructed
by Chernousov and Merkurjev (see \cite[Example~9.4]{CM04})
and is given by a $G$-homogeneous variety,
where $G$ is a product of two simple groups.
Our example
is given by a $G$-variety, where $G$ is a \emph{simple} group.
Apart from this, our example involves
an  (indecomposable) object in the category of integral Chow motives which is isomorphic over $\zz[\tfrac{1}{2}]$ and  $\zz[\tfrac{1}{3}]$ to the motive of the Severi-Brauer variety of a central simple division algebra of degree $5$, but not integrally.

The paper is organized as follows.
In section~\ref{stat} we state the main results. We 
then provide some technical facts
that are extensively used in the proofs
(section~\ref{notat}).
In the other sections we give
proofs of the results for varieties of type $\A_n$
(section~\ref{An}), of types $\B_n$ and $\C_n$ (section~\ref{Bn}), and exceptional
varieties of types $\G$ and $\FF$ (section~\ref{except}).
Section~\ref{krull} is devoted to the motivic decomposition
of generalized Severi-Brauer varieties.

\paragraph{Notation and Conventions}
By $G$ we denote an adjoint simple algebraic group of inner type
over $k$ and by $n$ its rank.
All varieties that appear in the paper are
projective $G$-homogeneous varieties over $k$. They can be considered as twisted forms of the varieties $G/P$, where $P$ is a parabolic subgroup of $G$. The Chow motive of a variety $X$ is denoted by $\M(X)$.
By $A$ we denote a central simple algebra over $k$ of index $\ind(A)$ and by $\SB(A)$ the corresponding Severi-Brauer variety. $I$ is always a right ideal of $A$ and $\rdim I$ stands for its reduced dimension. $V$ is a vector space over $k$.

By $\lambda=(\lambda_1,\lambda_2,\ldots,\lambda_l)$ we denote
a partition $\lambda_1\ge\lambda_2\ge\ldots\ge\lambda_l\ge 0$ with
$|\lambda|=\lambda_1+\lambda_2+\ldots+\lambda_l$.
Integers $d_1$, $d_2$, $\ldots$, $d_k$ always
satisfy the condition $1\le d_1<d_2<\ldots<d_k\le n$ and are
the dimensions of some flag. For each $i=0,\ldots,k$
we define $\delta_i$ to be the difference $d_{i+1}-d_i$ (assuming here
$d_0=0$ and $d_{k+1}=n+1$).

\section{Statements of Results}\label{stat}

We follow \cite[Appendix]{MPW96} and \cite{CG05}
for the description
of projective $G$-homogeneous varieties that appear below.
According to the type of the group $G$, we obtain the following results.


\paragraph{\large $\A_n$:}
In this case $G=\PGL_1(A)$, where
$A$ is a central simple algebra of degree $n+1$, $n>0$,
and the set of points of a projective $G$-homogeneous variety $X$ can be identified
with the set of flags of (right) ideals
$$
X(d_1,\ldots,d_k)=\{I_{1}\subset I_{2}\subset \ldots\subset I_{k} \subset A\}$$
of fixed reduced dimensions $1\le d_1<d_2<\ldots<d_k\le n$. Observe that this variety is a twisted form of $G/P$, where $P$ is the standard parabolic subgroup corresponding to the simple roots on the Dynkin diagram, numbered by $d_i$.
$$
\xymatrix@M=0pt{
\circ \ar@{}[d] |-{1} \ar@{-}[r]&
\circ \ar@{}[d] |-{2} \ar@{-}[r]&
\circ \ar@{}[d] |-{3} \ar@{.}[r]&
\circ \ar@{}[d] |-{n-2} \ar@{-}[r]&
\circ\ar@{}[d] |-{n-1}\ar@{-}[r]&
\circ\ar@{}[d] |-{n} \\
 & & & & &
}
$$

The following result reduces the computation of the motive of $X$
to the motives of ``smaller'' flags.

\begin{thm}\label{mainA}
Suppose that
$\gcd(\ind(A),d_1,\ldots,\hat{d}_m,\ldots,d_k)=1$, then
$$
\M(X(d_1,\ldots,d_k))\simeq \bigoplus_{\lambda}
\M(X(d_1,\ldots,\hat{d}_m,\ldots,d_k))(\delta_m\delta_{m-1}-|\lambda|),
$$
where the sum is taken over all partitions
$\lambda=(\lambda_1,\ldots,\lambda_{\delta_{m-1}})$ such that \\
$\delta_m\ge\lambda_1\ge\ldots\ge\lambda_{\delta_{m-1}}\ge 0$.
\end{thm}

\begin{proof} See \ref{proofmainA}.\end{proof}

As a consequence, for the variety of complete flags we obtain

\begin{cor}\label{corBorel}
The motive of the variety $X=X(1,\ldots,n)$ of complete
flags is isomorphic to
$$
\M(X)\simeq \bigoplus_{i=0}^{n(n-1)/2} \M(\SB(A))(i)^{\oplus a_i},
$$
where the
$a_i$ are the coefficients of the polynomial
$\varphi_n(z)=\sum_i a_iz^i=\prod_{k=2}^n \tfrac{z^k-1}{z-1}$.
\end{cor}
\begin{proof}
Apply Theorem~\ref{mainA} recursively to the sequence of varieties
$X(1,\ldots,n)$, $X(1,\ldots,n-1)$, $\ldots$, $X(1,2)$ and $X(1)=\SB(A)$.
\end{proof}

Another interesting example is the ``incidence'' variety $X(1,n)$:
\begin{cor}
The motive of $X(1,n)$ is isomorphic to
$$
\M(X(1,n))\simeq\bigoplus_{i=0}^{n-1}\M(\SB(A))(i).
$$
\end{cor}

In order to complete the picture we need to know
how to decompose the motive of a ``minimal'' flag,
i.e., a generalized Severi-Brauer variety.

Note that for some rings of coefficients (fields,
discrete valuation rings)
one easily obtains
the desired decomposition by using Krull-Schmidt Theorem
(the uniqueness of a direct sum decomposition).
More precisely, consider the subcategory $\M(G,R)$
of the category of motives with coefficients in a ring $R$
that is the pseudo-abelian completion of the category
of motives of projective $G$-homogeneous varieties
(see \cite[section~8]{CM04}).
Then we have the following

\begin{prop}\label{genSB}
Let $X(d)=\SB_d(A)$, $1< d< n$, be a generalized
Severi-Brauer variety for a central simple algebra $A$ of degree $n+1$
such that $\gcd(\ind(A),d)=1$.
Let $R$ be a ring such that
Krull-Schmidt Theorem holds in the category $\M(G,R)$.
Then the motive of $\SB_d(A)$ with
coefficients in $R$
is isomorphic to
$$
\M(\SB_d(A))\simeq \bigoplus_{i\in \mathcal{I}}
\M(\SB(A))(i)^{\oplus a_i},
$$
where the integers $a_i$ are the coefficients of the polynomial
$\tfrac{\varphi_{n}(z)}{\varphi_d(z)\varphi_{n+1-d}(z)}$ at terms $z^i$ and
the set of indices $\mathcal{I}=\{i\mid a_i\ne 0\}$.
\end{prop}

\begin{proof} See \ref{proofgenSB}.\end{proof}

It turns out  that the motives of some
generalized Severi-Brauer varieties with integral coefficients can still be
decomposed.

\begin{thm}\label{mainZ}
Let $\SB_2(A)$ be a generalized Severi-Brauer variety for
a division algebra of degree $5$.
Then there is an isomorphism
$$
\M(\SB_2(A))\simeq \Sb\oplus \Sb(2),
$$
where the motive $\Sb$ i) taken with $R$-coefficients, where $2$ or $3$ is invertible in a ring $R$, becomes isomorphic
to the motive of $\SB(A)$; \; ii) is indecomposable over $\zz$; \; iii)
is not isomorphic to the motive of $\SB(A)$ over $\zz$.
\end{thm}

\begin{proof} See \ref{proofmainZ}. \end{proof}

\begin{rem} It is expected that the case i) of the Theorem can also be proven
using techniques dealing with Norm varieties. Namely, if $A$ is cyclic, then
it corresponds to a symbol in $K^M_2(k)/5$ which splits the variety $\SB_2(A)$ 
(see \cite{RS}). By the results of Voevodsky \cite{Vo03} the motive of $\SB_2(A)$
splits into direct summands. 
\end{rem}

As an immediate consequence of Theorems~\ref{mainA} and \ref{mainZ}
we obtain

\begin{cor}\label{cex}
The Krull-Schmidt Theorem fails in the category of motives
$\M(\PGL_1(A),\zz)$ where $A$ is a division algebra of degree $5$.
\end{cor}

\begin{proof} Apply Theorem~\ref{mainA} recursively
to the sequences of varieties $X(1,2),X(1)$ and $X(1,2),X(2)$, where
$X(1,2)$ is the twisted flag $G$-variety for $G=\PGL_1(A)$.
We obtain two decompositions of the motive of $X(1,2)$
$$
\bigoplus_{i=0}^3 \M(\SB(A))(i)\simeq\M(X(1,2))\simeq
\M(\SB_2(A))\oplus \M(\SB_2(A))(1).
$$

Apply now Theorem~\ref{mainZ} to the components
of the second decomposition.
We obtain two different decompositions of the motive $\M(X(1,2))$
into indecomposable objects
$$
\bigoplus_{i=0}^3 \M(\SB(A))(i)\simeq\M(X(1,2))\simeq
\bigoplus_{i=0}^3 \mathcal{H}(i).
$$
This finishes the proof of the corollary.
\end{proof}

\begin{rem} Observe that the counter-example provided by Chernousov and Merkurjev
(see \cite[Ex.9.4]{CM04}) is the product of two Severi-Brauer
varieties $X=\SB(A)\times\SB(B)$ which is a $G$-homogeneous variety
for the \emph{semi-simple} group $G=\PGL_1(A)\times\PGL_1(B)$, where
$A$ and $B$ are two division algebras satisfying some conditions.
The example that we provide, i.e., the flag $X(1,2)$,
is a $G$-homogeneous variety for the \emph{simple} group $G=\PGL_1(A)$.
\end{rem}


\paragraph{\large $\B_n$:}

We assume that the characteristic of the base
field $k$ is different from $2$.
It is known that $G={\mathrm O}^+(V,q)$,
where $(V,q)$ is a regular quadratic space of dimension $2n+1$, $n>0$, and
projective $G$-homogeneous varieties can be described as flags
of totally $q$-isotropic subspaces
$$
X(d_1,\ldots,d_k)=\{V_1\subset \ldots\subset V_k\subset V\}.
$$
of fixed dimensions $1\le d_1< \ldots <d_k\le n$. Observe that this variety is a twisted form of $G/P$, where $P$ is the standard parabolic subgroup corresponding to the simple roots on the Dynkin diagram, numbered by $d_i$.
$$
\xymatrix@M=0pt{
\circ \ar@{}[d] |-{1} \ar@{-}[r]&
\circ \ar@{}[d] |-{2} \ar@{-}[r]&
\circ \ar@{}[d] |-{3} \ar@{.}[r]&
\circ \ar@{}[d] |-{n-2} \ar@{-}[r]&
\circ\ar@{}[d] |-{n-1}\ar@{=}[r] |-{>}&
\circ\ar@{}[d] |-{n} \\
 & & & & &
}
$$
The following result shows that some motives of flag
varieties can be decomposed into a direct sum
of twisted motives of ``smaller'' flags.

\begin{thm}\label{mainB}
Suppose that $m<k$, then
$$
\M(X(d_1,\ldots,d_k))\simeq\bigoplus_{\lambda}
\M(X(d_1,\ldots,\hat{d}_m,\ldots,d_k))(\delta_m\delta_{m-1}-|\lambda|),
$$
where the sum is taken over all partitions
$\lambda=(\lambda_1,\ldots,\lambda_{\delta_{m-1}})$ such that
$\delta_m\ge\lambda_1\ge\ldots\ge\lambda_{\delta_{m-1}}\ge 0$.
\end{thm}

\begin{proof} See \ref{proofmainB}. \end{proof}

In particular, for the variety of complete flags we obtain a formula similar to the one of Corollary 
\ref{corBorel}.

\begin{cor}
The motive of the variety of complete flags
$X=X(1,2,\ldots,n)$ is isomorphic to
$$
\M(X)\simeq \bigoplus_{i=0}^{n(n-1)/2} \M(X(n))(i)^{\oplus a_i},
$$
where the
$a_i$ are the coefficients of the polynomial
$\varphi_n(z)=\sum_i a_iz^i=\prod_{k=2}^n \tfrac{z^k-1}{z-1}$.
and $X(n)$ is the twisted form of the maximal orthogonal Grassmannian.
\end{cor}


\paragraph{\large $\C_n$:}

We assume that the characteristic of the base
field $k$ is different from $2$.
In this case $G=\Aut(A,\sigma)$,
where $A$ is a central simple algebra of degree $2n$, $n\ge2$,
with an involution $\sigma$ of symplectic type on $A$, and
a projective $G$-homogeneous variety
can be described as the set of flags of (right) ideals
$$
X(d_1,\ldots,d_k)=\{I_1\subset\ldots\subset I_k\subset A
\mid I_i\subseteq I_i^\bot\}
$$
of fixed reduced dimensions $1\le d_1<\ldots<d_k\le n$,
where $I^\bot=\{x\in A\mid\sigma(x)I=0\}$
is the right ideal of reduced dimension $2n-\rdim I$. Observe that this variety is a twisted form of $G/P$, where $P$ is the standard parabolic subgroup corresponding to the simple roots on the Dynkin diagram, numbered by $d_i$.
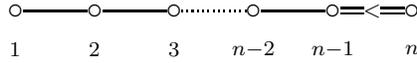
$$
\xymatrix@M=0pt{
\circ \ar@{}[d] |-{1} \ar@{-}[r]&
\circ \ar@{}[d] |-{2} \ar@{-}[r]&
\circ \ar@{}[d] |-{3} \ar@{.}[r]&
\circ \ar@{}[d] |-{n-2} \ar@{-}[r]&
\circ\ar@{}[d] |-{n-1}\ar@{=}[r] |-{<}&
\circ\ar@{}[d] |-{n} \\
 & & & & &
}
$$
Again, the motives of some flag
varieties can be decomposed into a direct sum
of twisted motives of ``smaller'' flags.

\begin{thm}\label{mainC}
Suppose that $d_i$ is odd for some $i<k$ and $d_k-d_{k-1}=1$.
Then
$$
\M(X(d_1,\ldots,d_k))\simeq
\bigoplus_{i=0}^{2n-2d_{k-1}-1}\M(X(d_1,\ldots,d_{k-1}))(i).
$$
\end{thm}

In particular, for the variety of complete flags we obtain

\begin{cor}
The motive of the variety of complete flags
$X=X(1,2,\ldots,n)$ is isomorphic to
$$
\M(X(1,\ldots,n))\simeq
\bigoplus_{i=0}^{n(n-1)} \M(\SB(A))(i)^{\oplus a_i},
$$
where $a_i$ are the coefficients of the polynomial
$\psi_n(z)=\prod_{k=1}^{n-1}\tfrac{z^{2k}-1}{z-1}$.
\end{cor}

\paragraph{\large $\G$:}
We suppose that the characteristic of $k$ is not $2$.
It is known that $G=\Aut(C)$,
where $C$ is a Cayley algebra over $k$.
By an \emph{$i$-space}, where $i=1,2$,
we mean an $i$-dimensional subspace $V_i$ of $C$
such that $uv=0$ for every $u,v\in V_i$.
The only flag variety corresponding to a non-maximal parabolic
is the variety of complete flags $X(1,2)$
which is described as follows
$$
X(1,2)=\{V_1\subset V_2\mid V_i\text{ is a $i$-subspace of }C\}.
$$
We enumerate the simple roots on the Dynkin diagram as follows:
$$
\xymatrix@M=0pt{
\circ\ar@{}[d] |-{1}\ar@3{-}[r] |-{<}& \circ\ar@{}[d] |-{2} \\
  &
}
$$

In this case we obtain
\begin{thm}\label{mainG} The motive of the variety of complete flags
$X=X(1,2)$ is isomorphic to
$$
\M(X)\simeq \M(X(2))\oplus\M(X(2))(1).
$$
\end{thm}
\begin{proof}
See \ref{proofmainG}
\end{proof}

Observe that by the result of Bonnet \cite{Bo03} the motives
of $X(1)$ and $X(2)$ are isomorphic (here $X(1)$ is a $5$-dimensional
quadric).

\paragraph{\large $\FF$:}
We suppose that the characteristic of $k$ is not $2$ and $3$.
It is known that $G=\Aut(J)$,
where $J$ is an exceptional Jordan algebra of dimension $27$ over $k$.
Set $\mathcal{I}=\{1,2,3,6\}$.
By an \emph{$i$-space}, $i\in \mathcal{I}$, we mean an $i$-dimensional subspace $V$ of $J$
such that every $u,v\in V$ satisfy the following condition:
$$
\tr(u)=0,\ u\times v=0,
\text{ and if }i<6\text{ then }u(va)=v(ua)\text{ for all }a\in J.
$$
A projective $G$-homogeneous variety
can be described as the set of flags of subspaces
$$
X(d_1,\ldots,d_k)=\{V_1\subset\ldots\subset V_k\mid V_i
\text{ is a $d_i$-subspace of }J\}.
$$
where the integers $d_1<\ldots<d_k$ are taken from the set $\mathcal I$.
Observe that this variety is a twisted form of $G/P$, where $P$ is the standard parabolic subgroup corresponding to the simple roots on the Dynkin diagram, numbered by $d_i$.
$$
\xymatrix@M=0pt{
\circ \ar@{}[d] |-{1} \ar@{-}[r]& \circ\ar@{}[d] |-{2}\ar@{=}[r] |-{<}& \circ\ar@{}[d] |-{3} \ar@{-}[r]& \circ\ar@{}[d] |-{6} \\
 & & &
}
$$
In this case we obtain

\begin{thm}\label{mainF}
Suppose that $m<k$ and either $d_{m+1}<6$ or $d_m=1$, then
$$
\M(X(d_1,\ldots,d_k))\simeq
\bigoplus_{\lambda}
\M(X(d_1,\ldots,\Hat d_m,\ldots,d_k))(\delta_m\delta_{m-1}-|\lambda|),
$$
where the sum is taken over all partitions
$\lambda=(\lambda_1,\ldots,\lambda_{\delta_{m-1}})$ such that
$\delta_m\ge\lambda_1\ge\ldots\ge\lambda_{\delta_{m-1}}\ge 0$.
\end{thm}

\begin{proof} See \ref{proofmainF}. \end{proof}


\section{Preliminaries}\label{notat}

In the present section, we introduce the category of Chow motives
following \cite{Ma68}. We formulate and prove the Grassmann Bundle Theorem.
At the end we recall the notion of functor of points following
\cite[section~8]{Ka01} and provide some examples.

\begin{ntt}[Chow motives]
Let $k$ be a field and $\Var_k$ be the category
of smooth projective varieties over $k$.
We define the category $\Cor_k$ of \emph{correspondences} over $k$.
Its objects are smooth projective varieties over $k$.
As morphisms, called correspondences,
we set $\Mor(X,Y):=\CH^{\dim X}(X\times Y)$.
For any two correspondences $\alpha\in \CH(X\times Y)$ and
$\beta\in \CH(Y\times Z)$ we define the composition
$\beta\circ\alpha\in \CH(X\times Z)$
$$
\beta\circ\alpha ={\pr_{13}}_*(\pr_{12}^*(\alpha)\cdot \pr_{23}^*(\beta)),
$$
where $\pr_{ij}$ denotes the projection
on the $i$-th and $j$-th factors of $X\times Y\times Z$
respectively and ${\pr_{ij}}_*$, ${\pr_{ij}^*}$ denote
the induced push-forwards and pull-backs for Chow groups.
Observe that the composition $\circ$
induces a ring structure on the abelian group $\CH^{\dim X}(X\times X)$.
The unit element of this ring is the class of the diagonal $\Delta_X$.

The pseudo-abelian completion of $\Cor_k$ is called the category
of \emph{Chow motives} and is denoted by $\M_k$.
The objects of $\M_k$
are pairs $(X,p)$, where $X$ is a smooth projective variety
and $p$ is a projector, that is, $p\circ p=p$.
The motive $(X,\Delta_X)$ will be denoted by $\M(X)$.
\end{ntt}

\begin{ntt}
By the construction $\M_k$ is a self-dual tensor additive category,
where the duality is given by the transposition of cycles $\alpha\mapsto \alpha^t$
and the tensor product is given by the usual fiber product
$(X,p)\otimes(Y,q)=(X\times Y, p\times q)$.
Moreover, the Chow functor
$\CH: \Var_k \to \zab$ (to the category of $\zz$-graded abelian groups)
factors through $\M_k$, i.e., one has the commutative diagram of functors
$$
\xymatrix{
\Var_k \ar[rr]^{\CH}\ar[rd]_{\Gamma}& & \zab \\
 & \M_k \ar[ru]_{R}&
}
$$
where $\Gamma: f\mapsto \Gamma_f$ is the graph and
$R:(X,p) \mapsto \im(p^*)$ is the realization.
\end{ntt}

\begin{ntt} Consider the morphism
$(\id,e):\pl^1\times\{pt\}\to \pl^1\times\pl^1$. Its image by means of
the induced push-forward $(\id,e)_*$ doesn't depend on the choice of a point
$e:\{pt\}\to \pl^1$
and defines the projector in $\CH^1(\pl^1\times\pl^1)$ denoted by $p_1$.
The motive $L=(\pl^1,p_1)$ is called the Lefschetz motive.
For a motive $M$ and an nonnegative integer $i$ we denote by $M(i)=M\otimes L^{\otimes i}$ its twist.
Observe that
$$
\Mor((X,p)(i),(Y,q)(j))=q\circ \CH^{\dim X+i-j}(X,Y)\circ p.
$$
\end{ntt}

We will extensively use the following fact
that easily follows from Manin's Identity
Principle and Grassmann Bundle Theorem for Chow groups \cite{Ful}.

\begin{prop}[Grassmann Bundle Theorem]\label{GrBundl}
Let $X$ be a variety over $k$ and
$\E$ be a vector bundle over $X$ of rank $n$.
Then the motive of the Grassmann bundle
$\Gr(d,\E)$ over $X$ is isomorphic to
$$
\M(\Gr(d,\E))\simeq \bigoplus_\lambda\M(X)(d(n-d)-|\lambda|),
$$
where the sum is taken over all partitions
$\lambda=(\lambda_1,\ldots,\lambda_d)$ such that
$n-d\ge\lambda_1\ge\ldots\ge\lambda_d\ge 0$.
\end{prop}

\begin{proof}
We follow notations of \cite{Ma68}. Denote $\Gr(d,\E)$ by $Y$ and the canonical projection of $Y$ to $X$ by $\psi$. It is known (see \cite[Prop.~14.6.5]{Ful}) that $\CH^*(Y)$ as $\CH^*(X)$-module (via $\psi^*$) has a basis consisting of elements $\De_\lambda\in\CH^{|\lambda|}(Y)$ parametrized by partitions $\lambda$. For any partition $\lambda$ denote by $\lambda^\op $ the partition defined by
$\lambda^\op _i=n-d-\lambda_{d+1-i}$. Then for every two partitions $\lambda$ and $\mu$ with $|\lambda|+|\mu|\le d(n-d)$ and any element $a\in\CH^*(X)$ the following duality formula holds (see \cite[Prop.~14.6.3]{Ful}):
$$
\psi_*(\De_\lambda\De_\mu\psi^*(a))=
\begin{cases}
&a,\ \mu=\lambda^\op ,\\
&0,\text{ otherwise.}
\end{cases}
$$
For every partition $\lambda$, we define an element $f_\lambda$ of $\End(\M(Y))$ inductively. For the unique partition $\lambda_{\max}$ with $|\lambda_{\max}|=d(n-d)$ set $f_{\lambda_{\max}}=c(\psi)\circ c(\psi)^t$. Now, by decreasing induction on $|\lambda|$, set
$$
f_\lambda=f_{\lambda_{\max}}\circ c_{\De_{\lambda^\op }}\circ(\Delta_Y-\sum_{|\mu|>|\lambda|}c_{\De_\mu}\circ f_\mu).
$$
Finally, set $p_\lambda=c_{\De_\lambda}\circ f_\lambda$.

Now, the duality formula implies that
$$
(p_{\lambda})_e(\sum_\mu\De_\mu\psi^*(a_\mu))=\De_\lambda\psi^*(a_\lambda).
$$
Therefore, by Manin's identity principle, the $p_\lambda$ form a complete orthogonal system of projectors.
Identification of their images with twisted motives of $X$ can be done as in \cite[\S 7]{Ma68} and we omit it.
\end{proof}

\begin{ntt}[Functors of Points]\label{functpoints}
In sections \ref{An}, \ref{Bn} and \ref{except}
we use the functorial language,
that is consider $k$-schemes as functors from the category of $k$-algebras
to the category of sets. Fix a scheme $X$.
By an $X$-algebra we mean a pair $(R,x)$,
where $R$ is a $k$-algebra and $x$ is an element of $X(R)$.
$X$-algebras form a category with obvious morphisms.
Morphisms $\varphi\colon Y\to X$ can be considered as functors
from the category of $X$-algebras to the category of sets, by sending a pair $(R,x)$ to its preimage in $Y(R)$.
\end{ntt}

\begin{ntt}  Let $X$ be a variety over $k$.
To any vector bundle $\mathcal{F}$ over $X$ we can associate
the Grassmann bundle $Y=\Gr(d,\mathcal{F})$. Fix an $X$-algebra $(R,x)$.
The value of the
functor corresponding to $\Gr(d,\mathcal{F})$ at $(R,x)$
is the set of direct summands of rank $d$ of the projective $R$-module
$\mathcal{F}_x\otimes_k R$ (cf. \cite[section~9]{Ka01}), where $\mathcal{F}_x=\F(R,x)$.
\end{ntt}

\section{Groups of type $\A_n$}\label{An}

The goal of the present section is to prove Theorem~\ref{mainA} and Proposition~\ref{genSB}.
We use the notation of section~\ref{stat}.

\begin{ntt}
Let $G$ be an adjoint group of inner type $\A_n$ defined over a field $k$.
It is well known that $G=\PGL_1(A)$, where $A$ is a central simple algebra of degree $n+1$ and
points of projective $G$-homogeneous varieties are flags of (right) ideals of $A$
$$
X(d_1,\ldots,d_k)=\{I_1\subset\ldots\subset I_k\subset A \mid \rdim I_i=d_i\}.
$$
\end{ntt}

For convenience we set $d_0=0$, $d_{k+1}=n+1$, $I_0=0$, $I_{k+1}=A$.

\begin{ntt}\label{condid}
The value of the functor of points corresponding to the variety $X(d_1,\ldots,d_k)$
at a $k$-algebra $R$ (see \ref{functpoints}) equals the set of all flags
$I_1\subset\ldots\subset I_k$ of right ideals of $A_R=A\otimes_k R$
having the following properties (see \cite[section~4]{IK00})
\begin{itemize}
\item the injection of $A_R$-modules $I_i\hookrightarrow A_R$ splits;
\item $\rdim I_i=d_i$.
\end{itemize}
\end{ntt}

\begin{ntt}
On the scheme $X=X(d_1,\ldots,d_k)$
there are ``tautological'' vector bundles $\J_i$, $i=0,\ldots,k+1$, of ranks
$(n+1)d_i$. The value of $\J_i$ on an $X$-algebra $(R,x)$,
where $x=(I_1,\ldots,I_k)$, is the ideal $I_i$ considered as a projective $R$-module.
The bundle $\J_i$ also has a structure of right $A_X$-module,
where $A_X$ is the constant sheaf of algebras on $X$ determined by $A$.

For every $m\in\{1,\ldots,k\}$ there exists an obvious morphism
\begin{align*}
X(d_1,\ldots,d_k)& \to X(d_1,\ldots,\Hat d_m,\ldots,d_k)\\
(I_1,\ldots,I_k)&\mapsto (I_1,\ldots,\Hat I_m,\ldots,I_k)
\end{align*}
that turns $X(d_1,\ldots,d_k)$ into a $X(d_1,\ldots,\Hat d_m,\ldots,d_k)$-scheme.
\end{ntt}

\begin{lem}\label{thmA}
Denote $X(d_1,\ldots,d_k)$ by $Y$ and $X(d_1,\ldots,\Hat d_m,\ldots,d_k)$ by $X$.
Assume there exists a vector bundle $\E$ over $X$ such that
$A_X\simeq\End_{\OO_X}(\E)$. Consider the vector bundle
$$
\F=\J_{m+1}\E/\J_{m-1}\E=\J_{m+1}/\J_{m-1}\otimes_{A_X}\E
$$
of rank $d_{m+1}-d_{m-1}$.
Then $Y$ as a scheme over $X$
can be identified with the Grassmann bundle $Z=\Gr(d_m-d_{m-1},\F)$ over $X$.
\end{lem}

\begin{proof}
We use essentially the same method as in \cite[Proposition 4.3]{IK00}.

Fix an $X$-algebra $(R,x)$ where $x=(I_1,\ldots,\Hat I_m,\ldots,I_k)$.
The fiber of $Y$ over $x$, i.e., the value at $(R,x)$,
can be identified with the set of all
ideals $I_m$ satisfying the conditions \ref{condid}
such that $I_{m-1}\subset I_m \subset I_{m+1}$.
The fiber of $Z$ over $x$ is the set of all $R$-submodules $N$ of
$\F_y=\F(R,y)$ such that the injection $N\hookrightarrow\F_y$ splits and
$\rk_R N=d_m-d_{m-1}$.

We define a natural bijection between
the fibers of $Y$ and $Z$ over $x$ as follows.

Consider the following mutually inverse bijections between the
set of all right ideals of reduced dimension $r$ in $A_R$
(satisfying \ref{condid}) and
the set of all direct summands of rank $r$ of
the $R$-module $\E_x$
\begin{align*}
\Phi:I&\mapsto I\E_x\\
\Psi:N&\mapsto\Hom_R(\E_x,N)\subset\End_R(\E_x)\simeq A_R
\end{align*}
Observe that these bijections preserve the respective inclusions of ideals and modules.
So ideals of reduced dimension $d_m$ between $I_{m-1}$
and $I_m$ correspond to
submodules of rank $d_m$ between $I_{m-1}\E_x$ and $I_{m+1}\E_x$,
and, therefore, to submodules of rank $d_{m+1}-d_{m-1}$ in
$I_{m+1}\E_x/I_{m-1}\E_x=\F_x$.
This gives the desired natural bijection on the fibers.
\end{proof}

\begin{lem}\label{propBundle}
Suppose that $\gcd(\ind(A),d_1,\ldots,d_k)=1$.
Then there exists a vector bundle $\E$ over $X=X(d_1,\ldots,d_k)$
of rank $n+1$
such that $A_X\simeq\End_{\OO_X}(\E)$.
\end{lem}

\begin{proof}
We have to prove that the class $[A_X]$ in $\Br(X)$ is trivial.
Since $X$ is a regular Noetherian scheme the canonical map
$$
\Br(X)\to\Br(K)
$$
where $K=k(X)$ is the function field of $X$, is injective by \cite[1.10]{Gr}
and \cite[Theorem~7.2]{AG60}.
So it is enough to prove that
$A\otimes_k K$ splits.
But the generic point of $X$ defines a flag of ideals of $A\otimes_k K$
of reduced dimensions $d_1,\ldots,d_k$. Since
the index $\ind(A\otimes_k K)$ divides $d_1,\ldots, d_k$ and $\ind A$,
by the assumption of the lemma it must be equal to $1$.
So $A\otimes_k K$ is split and
this finishes the proof of the lemma.
\end{proof}

\begin{rem}
In the case $d_1=1$ one can take $\E=\J_1^\vee$.
\end{rem}

\begin{rem}
It can be shown using the Index Reduction Formula (see \cite{MPW96}) that
the condition on the $\gcd$ is necessary and sufficient for the central simple algebra
$A_{k(X)}$ to be split.
\end{rem}

We are now ready to finish the proof of \ref{mainA}.

\begin{ntt}[Proof of Theorem \ref{mainA}]\label{proofmainA}
By Lemma~\ref{propBundle} there exists a vector bundle $\E$ over
$X=X(d_1,\ldots,\hat{d}_m,\ldots,d_k)$ of rank $n+1$ such that
$A_X\simeq\End_{\OO_X}(\E)$.
By Lemma~\ref{thmA} we conclude that $Y=X(d_1,\ldots,d_k)$ is a Grassmann bundle
over $X$. Now by Proposition~\ref{GrBundl} we obtain the isomorphism of \ref{mainA}.
\end{ntt}

\begin{rem}
Note that the assumption of Theorem~\ref{mainA} on the reduced dimensions $d_1,\ldots,d_k$
is essential.
Indeed, suppose the Theorem holds for any twisted flag variety.
Consider the flag $X=X(1,d)$ with $\gcd(\ind(A),d)>1$.
Then we have an isomorphism of motives
$$
\M(X)\simeq\bigoplus_{i=0}^{d-1}\M(\SB_d(A))(i)
$$
which appears after applying Theorem to the flags $X(1,d)$ and $X(d)$.
Consider the group $\CH_0(X)=\Mor_{\M_k}(\M(pt),\M(X))$.
The isomorphism above induces the isomorphism of groups
\begin{align*}
\Coker(\CH_0(X)\xra{\res}\CH_0(X_{k_s}))&\cong \Coker(\CH_0(\SB_d(A))\xra{\res}\CH_0(\Gr(d,n+1)))\\
&\cong \zz/(\tfrac{\ind(A)}{\gcd(\ind(A),d)})\zz,
\end{align*}
where $\res$ is the pull-back induced by the scalar extension $k_s/k$
(here $k_s$ denotes the separable closure of $k$) and the last isomorphism
follows by \cite[Theorem~3]{Bl91}.
On the other hand, applying Theorem~\ref{mainA} to the flags $X(1,d)$ and $X(1)$ we
obtain the isomorphism
$$
\M(X)\simeq\bigoplus_\lambda\M(\SB(A))((n+1-d)(d-1)-|\lambda|)
$$
which induces the isomorphism of groups
\begin{align*}
\Coker(\CH_0(X)\xra{\res}\CH_0(X_{k_s}))& \cong\Coker(\CH_0(\SB(A))\xra{\res}\CH_0(\pl^n)) \\
& \cong \zz/\ind(A)\zz,
\end{align*}
that leads to a contradiction.
\end{rem}

We now prove Proposition~\ref{genSB}.

\begin{ntt}[Proof of Proposition~\ref{genSB}]\label{proofgenSB}
Let $G=\PGL_1(A)$ and let $\M(G,R)$ be the symmetric tensor category
of motives of $G$-homogeneous varieties
with coefficients in a ring $R$ for which the Krull-Schmidt theorem holds.
It is the case, e.g., when $R$ is a field or, more general, a discrete valuation ring
(see \cite[Theorem~9.6]{CM04}).

Consider the $G$-homogeneous variety $X(1,d)$, $1<d<n$.
Apply Theorem~\ref{mainA} to the sequences of flags $X(1,d)$, $X(d)$ and $X(1,d)$, $X(1)$.
We obtain two isomorphisms in $\M(G,R)$
\begin{equation}\tag{*}
\bigoplus_{i=0}^{d-1}\M(\SB_d(A))(i)\simeq\M(X)\simeq
\bigoplus_\lambda\M(\SB(A))((n+1-d)(d-1)-|\lambda|),
\end{equation}
where the sum on the right hand side is taken over all partitions
$\lambda=(\lambda_1,\ldots,\lambda_{d-1})$ such that $n+1-d\ge \lambda_1\ge \ldots\ge \lambda_{d-1}\ge 0$.
Since Krull-Schmidt Theorem holds in $\M(G,R)$, the motive
$\SB(A)$ has a unique decomposition into the direct sum of indecomposable
objects $H_i$, $i\in\mathcal{I}$, and their twists
$$
\M(\SB(A))\simeq \bigoplus_{i\in\mathcal{I}} (\oplus_{j\in \mathcal{J}_i} H_i(j)).
$$

Consider the subcategory $\M(G,R)_\mathcal{I}$
additively generated by the
motives $H_i$, $i\in\mathcal{I}$, and their twists.
The abelian group of isomorphism classes of objects
of this category can be equipped with a structure
of a free module over the polynomial ring $R[z]$.
Namely, multiplication by $z$ is given by the twist.
Clearly, the classes $[H_i]$, $i\in\mathcal{I}$,
form the basis of this $R[z]$-module.

By (*) we have $\M(\SB_d(A))\in \M(G,R)_{\mathcal{I}}$ and
the isomorphisms (*) can be rewritten as
\begin{align*}
\frac{z^d-1}{z-1}[\SB_d(A)]&=
\frac{\varphi_{n}(z)}{\varphi_{d-1}(z)\varphi_{n+1-d}(z)}[\SB(A)]\\
&=\frac{z^d-1}{z-1}
\frac{\varphi_{n}(z)}{\varphi_{d}(z)\varphi_{n+1-d}(z)}[\SB(A)]
\end{align*}
where $\varphi_n(z)=\prod_{k=2}^n\tfrac{z^k-1}{z-1}$.
This immediately implies the equality
$$
[\SB_d(A)]=\frac{\varphi_{n}(z)}{\varphi_{d}(z)\varphi_{n+1-d}(z)}[\SB(A)],
$$
i.e., the isomorphism in $\M(G,R)_{\mathcal{I}}$
between $\M(\SB_d(A))$ and the respective sum of twists of $\M(\SB(A))$.
This finishes the proof of the proposition.
\end{ntt}


\section{Motivic decomposition of $\SB_2(A)$}\label{krull}

This section is devoted to the proof of Theorem~\ref{mainZ}.

First, we need to recall some properties of rational
cycles on projective homogeneous varieties.

\begin{ntt} 
Let $G$ be a split linear algebraic group over $k$.
Let $X$ be a projective $G$-homogeneous variety, i.e., 
$X=G/P$, where $P$ is a parabolic subgroup of $G$.
The abelian group structure of $\CH(X)$, as well as
its ring structure, is well-known.
Namely, $X$ has a cellular filtration and the generators of Chow groups of the bases of this filtration
correspond to the free additive generators of $\CH(X)$ (see \cite{Ka01}).
Note that the product of two projective homogeneous varieties
$X\times Y$ has a cellular filtration as well,
and $\CH^*(X\times Y)\cong \CH^*(X)\otimes \CH^*(Y)$
as graded rings.
The correspondence product of two cycles 
$\alpha=f_\alpha \times g_\alpha \in \CH(X\times Y)$ and
$\beta=f_\beta \times g_\beta \in \CH(Y\times X)$ is given
by (cf. \cite[Lem.~5]{Bo03})
$$
(f_\beta\times g_\beta)\circ(f_\alpha\times g_\alpha)=
\deg(g_\alpha \cdot f_\beta)(f_\alpha\times g_\beta),
$$
where $\deg: \CH(Y)\to \CH(\{pt\})=\zz$ is the degree map.
\end{ntt}

\begin{ntt}
Let $X$ be a projective variety of dimension $n$ over a field $k$.
Let $k_s$ be the separable closure of the field $k$.
Consider the scalar extension $X_s=X\times_k k_s$.
We say a cycle $J\in \CH(X_s)$ is {\it rational}
if it lies in the image of the pull-back homomorphism
$\CH(X)\to \CH(X_s)$.
For instance, there is an obvious rational cycle $\Delta_{X_s}$ on
$\CH^n(X_s\times X_s)$ that is given by the diagonal class.
Clearly, linear combinations,
intersections and correspondence products of rational cycles
are rational.
\end{ntt}

\begin{ntt}\label{exproj}
We will use the following fact (see \cite[Cor.~8.3]{CGM})
that follows from the Rost Nilpotence Theorem.
Let $p_s$ be a non-trivial rational projector in $\CH^n(X_s\times X_s)$,
i.e., $p_s\circ p_s=p_s$. Then there exists a non-trivial projector $p$
on $\CH^n(X\times X)$ such that $p\times_k k_s=p_s$.
Hence, the existence of a non-trivial rational projector $p_s$ on
$\CH^n(X_s\times X_s)$ gives rise to
the decomposition of the Chow motive of $X$
$$
\M(X)\cong(X,p)\oplus (X,\Delta_X-p)
$$
Our goal is to find such a projector
in the case $X=\SB_d(A)$.
\end{ntt}

\begin{ntt}\label{twistisom}
An isomorphism between $(X,p)(m)$ and $(Y,q)(l)$ can be given by correspondences $j_1\in\CH^{\dim X-l+m}(X\times Y)$ and $j_2\in\CH^{\dim Y-m+l}(Y\times X)$ such that
$q\circ j_1=j_1\circ p$, $p\circ j_2=j_2\circ q$, $j_1\circ j_2=q$, $j_2\circ j_1=p$. If $X$ and $Y$ lie in the category $\M(G,\zz)$ then by the Rost nilpotence theorem (see \cite[Theorem~8.2]{CM04} and \cite[Corollary~8.4]{CGM}) it suffices to give a rational $j_1$ and some $j_2$ satisfying these conditions over separable closure (note that $j_2$ will automatically be rational).
\end{ntt}

We now recall some properties of Grassmann varieties and
describe their Chow rings.

\begin{ntt}
Consider the Grassmann variety $\Gr(d,n+1)$, $1\le d \le n$,
of $d$-planes in the $(n+1)$-dimensional affine space.
It has the dimension $d(n+1-d)$. The twisted form of it is a generalized
Severi-Brauer variety $\SB_d(A)$, where $A$ is a central simple algebra
of degree $n+1$.
For any two integers $d$ and $d'$, $1\leq d,d'\leq n$,
there is the fiber product diagram
\begin{equation}\label{segre}
\xymatrix{
\Gr(d,n+1)\times\Gr(d',n+1)\ar[r]^-{\Seg}\ar[d]&\Gr(dd',(n+1)^2)
\ar[d]^{\simeq}\\
\SB_d(A)\times\SB_{d'}(A^\op )\ar[r]^-{\Seg}&\SB_{dd'}(A\otimes A^\op )
}
\end{equation}
where the horizontal arrows are Segre embeddings given
by the tensor product of ideals (resp. linear subspaces) and the vertical
arrows are canonical maps induced by the scalar extension $k_s/k$
(here $k_s$ is the separable closure of $k$).
Observe that the right vertical arrow is an isomorphism
since $A\otimes A^\op $ splits.
\end{ntt}

\begin{ntt}
The diagram (\ref{segre}) induces the commutative diagram of rings
\begin{equation}
\xymatrix{
\CH(\Gr(d,n+1) \times\Gr(d',n+1)) & \CH(\Gr(dd',(n+1)^2)) \ar[l]_-{\Seg^*}\\
\CH(\SB_d(A)\times\SB_{d'}(A^\op ))\ar[u]_{\res} & \CH(\SB_{dd'}(A\otimes A^\op ))
\ar[l]_-{\Seg^*}\ar[u]^{\simeq}_{\res}
}
\end{equation}
where all maps are the induced pull-backs.
Consider a vector bundle $E$ over $\Gr(dd',(n+1)^2)$.
It is easy to see that
the pull-back of the total Chern class $\Seg^*(c(E))$ is a rational cycle on
$\CH(\Gr(d,n+1)\times\Gr(d',n+1))=\CH(\Gr(d,n+1)) \otimes \CH(\Gr(d',n+1))$.
In particular, if $E=\tau_{dd'}$ is the
tautological bundle of $\Gr(dd',(n+1)^2)$
we obtain the following
\end{ntt}

\begin{lem}\label{ratc}
The total Chern class
$c(\pr_1^*\tau_d\otimes\pr_2^*\tau_{d'})$ of the tensor product
of the pull-backs (induced by the projection maps)
of the tautological bundles $\tau_d$ and $\tau_{d'}$ of $\Gr(d,n+1)$ and
$\Gr(d',n+1)$ respectively is rational.
\end{lem}

From now on we restrict ourselves to
the case $n=4$, $d=2$ and $d'=1$, i.e.,
to the Grassmannian $\Gr(2,5)$ and the projective space $\pl^4=\Gr(1,5)$.

\begin{ntt}
We describe the generators and relations of the Chow ring $\CH(\Gr(2,5))$
following \cite[section~14.7]{Ful}.
Set $\sigma_m=c_m(Q)$, $m=1,2,3$,
where $Q=\OO^5/\tau_2$ is the universal quotient
bundle of rank $3$ over $\Gr(2,5)$.
It is known that the elements $\sigma_m$
generate the Chow ring $\CH(\Gr(2,5))$. More precisely, as an abelian group
this ring is generated by the
Schubert cycles $\De_\lambda(\sigma)$
that are parametrized by all partitions $\lambda=(\lambda_1,\lambda_2)$
such that $3\ge\lambda_1\ge\lambda_2\ge 0$.
In particular, $\sigma_m=\De_{(m,0)}$, $m=1,2,3$.
For other generators we set the following notation
$g_2=\De_{(1,1)}$, $g_3=\De_{(2,1)}$, $h_4=\De_{(3,1)}$, $g_4=\De_{(2,2)}$,
$g_5=\De_{(3,2)}$, $pt=\De_{(3,3)}$.
These generators corresponds to the vertices of the Hasse diagram of $\Gr(2,5)$
$$
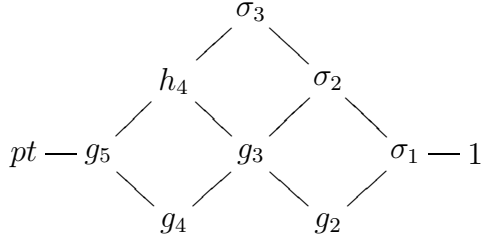
\xymatrix@=1em{
&&&\sigma_3\ar@{-}[ld]\ar@{-}[rd]\\
&&h_4\ar@{-}[ld]\ar@{-}[rd]&&\sigma_2\ar@{-}[ld]\ar@{-}[rd]\\
pt\ar@{-}[r]&g_5\ar@{-}[rd]&&g_3\ar@{-}[ld]\ar@{-}[rd]&&\sigma_1\ar@{-}[ld]\ar@{-}[r]&1\\
&&g_4&&g_2
}
$$
The multiplication rules can be determined using Pieri formulae
$$
\De_\lambda \cdot \sigma_m = \sum_{\mu} \De_\mu,
$$
where the sum is taken over all partitions $\mu=(\mu_1,\mu_2)$ such that
$3\ge \mu_1 \ge \lambda_1 \ge \mu_2 \ge \lambda_2\ge 0$.
\end{ntt}

\begin{ntt} Consider the tautological bundle $\tau_2$ of the
Grassmannian $\Gr(2,5)$. Its total Chern class is 
$$
c(\tau_2)=c(Q)^{-1}=
\frac{1}{1+\sigma_1+\sigma_2+\sigma_3}=1-\sigma_1-\sigma_2+\sigma_1^2+\ldots
$$
where the rest consists of the summands of degree greater than $2$.
Hence, we obtain $c_1(\tau_2)=-\sigma_1$
and $c_2(\tau_2)=-\sigma_2+\sigma_1^2=g_2$.
\end{ntt}

\begin{ntt} The Chow ring of the projective space $\pl^4$
can be identified with the factor ring $\zz[H]/(H^5)$,
where $H=c_1(\OO(1))$ is the class of a hyperplane section.
Thus, the first Chern class of the tautological bundle of $\pl^4$ equals to
$c_1(\tau_1)=c_1(\OO(-1))=-H$.
\end{ntt}

We are now ready to prove Theorem~\ref{mainZ}.

\begin{ntt}[Proof of Theorem~\ref{mainZ}]\label{proofmainZ}
By Lemma~\ref{ratc} we obtain
the following rational cycles in  $\CH^*(\Gr(2,5)\times\pl^4)$
\begin{align*}
r&=c_1(\pr_1^*(\tau_2)\otimes\pr_2^*(\tau_1))=c_1(\pr_1^*(\tau_2)+2c_1(\pr_2^*(\tau_1))
=-\sigma_1\times 1-2(1\times H),\\
\rho&=c_2(\pr_1^*(\tau_2)\otimes\pr_2^*(\tau_1))=c_2(\pr_1^*(\tau_2))+c_1(\pr_1^*(\tau_2))c_1(\pr_2^*(\tau_1))+c_1(\pr_2^*(\tau_1))^2\\
&=g_2\times 1+\sigma_1\times H+1\times H^2
\end{align*}
\end{ntt}

For two cycles $x$ and $y$ we shall write $x=_5y$ if there exists a cycle $z$ such that $x-y=5z$.
Note that $=_5$ is an equivalence relation
that preserves rationality of cycles.
Then the following cycles are rational
\begin{align*}
\rho^2&=_5 1\times H^4+2\sigma_1\times H^3+(\sigma_2+3g_2)\times H^2+2g_3\times H+g_4\times 1,\\
\rho^3&=_5(3\sigma_2+g_2)\times H^4+(\sigma_3+3g_3)\times H^3+(g_4+3h_4)\times H^2+3g_5\times H+pt\times 1.
\end{align*}
and a direct computation shows that
$\rho^3\circ(\rho^2)^t=_5\Delta_{\pl^4}$ is the class of the diagonal in
$\CH^4(\pl^4\times\pl^4)$.

Consider the composition
\begin{align*}
(\rho^2)^t\circ\rho^3=_5&(3\sigma_2+g_2)\times g_4+(2\sigma_3+g_3)\times g_3\\
&+(g_4+3h_4)\times (\sigma_2-2g_2)+g_5\times \sigma_1+pt\times 1.
\end{align*}
Note that the right-hand side is a rational projector (over $\zz$)
and, therefore, by Rost Nilpotence Theorem
(see \cite[Corollary~8.3]{CGM})
has a form $p\times_k k_s$ where $p$ is a projector
in $\End(\M(\SB_2(A)))$.
The latter determines an object
$(\SB_2(A),p)$ in the category of motives
(actually in $\M(G,\zz)$) which we denoted by $\Sb$.

Set $q=\Delta_{\SB_2(A)}-p$. We then show that
$$
(\M(\SB_2(A)),q)\simeq(\M(\SB_2(A)),p^t)\simeq\Sb(2),
$$
which gives the claimed decomposition $\M(\SB_2(A))\simeq\Sb\oplus\Sb(2)$.

Observe that an isomorphism $(\SB_2(A),q)\simeq(\SB_2(A),p^t)$
is given by the two mutually inverse motivic isomorphisms
$p^t_s\circ q_s$ and $q_s\circ p^t_s$
over $k_s$ which are rational.
An isomorphism $\Sb(2)\simeq(\M(\SB_2(A)),p^t)$ is given by the
following two cycles
\begin{align*}
j_1&=(3\sigma_2+g_2)\times pt-(2\sigma_3+g_3)\times g_5\\
&+(g_4+3h_4)\times (g_4+3h_4)-g_5\times (2\sigma_3+g_3)+
pt\times (3\sigma_2+g_2),\\
j_2&=1\times g_4-\sigma_1\times g_3+(\sigma_2-2g_2)\times (\sigma_2-2g_2)-g_3\times \sigma_1+
g_4\times 1.
\end{align*}
Note that $j_1$ is rational, since $j_1=_5 (1\times(3\sigma_2+g_2))p$, and
$1\times(3\sigma_2+g_2)=_5 3(\rho+r^2)^t\circ\rho^2$ is rational.

We now check the properties of the motive $\Sb$.

i) To prove that $\Sb$ is isomorphic to $\M(\SB(A))$ over $R$ it is enough to provide rational cycles $\alpha\in\CH^4(\pl^4\times \Gr(2,5))$ and $\beta\in\CH^6(\Gr(2,5)\times\pl^4)$ such that $\alpha\circ\beta=p_s$ and $\beta\circ\alpha=\Delta_{\pl^4}$. If $2$ is invertible in $R$, set
$$
\alpha=(\rho^2)^t,\quad\beta=\rho^3-\tfrac{5}{2}(g_5\times H+g_3\times H^3).
$$
If $3$ is invertible, set
$$
\alpha=(\rho^2)^t-\tfrac{5}{3}(H\times g_3+H^3\times\sigma_1)-5H^2\times g_2,
\quad\beta=\rho^3.
$$

ii) We show that $\Sb$ is indecomposable over $\zz$. 
It is enough to prove that it is indecomposable modulo $5$. 
Observe that by item i) the motive $\Sb$ with $\zz/5\zz$-coefficients
is isomorphic to $\M(\SB(A))$. Assume that $\M(\SB(A))$ is decomposable modulo 5. 
Let $\mu$ be a non-trivial projector in $\CH(\SB(A)\times\SB(A))/5$.
Then $\mu_s=\sum_i{\pm}H^i\times H^{5-i}$. Chose any inverse image $\nu$ of $\mu$ in $\CH(\SB(A)\times\SB(A))$. Then $\nu_s=_5\sum_i{\pm}H^i\times H^{5-i}$ is rational, and so we have a rational non-trivial projector over $\zz$, that contradicts \cite[theorem~2.2.1]{Ka95}.

iii) We claim that $\Sb\not\simeq\M(\SB(A))$ over $\zz$.
Indeed, assume there is a motivic isomorphism
$\varphi\colon\M(\SB(A))\xra{\simeq}\Sb$.
Since $\Sb$ is a direct summand of $\M(\SB_2(A))$,
there exist maps $\alpha'$, $\beta'$ such that $\alpha'\circ \beta'=p$,
$\beta'\circ\alpha'=\Delta_{\Sb}$ (see the diagram).
$$
\xymatrix{
{\M(\SB_2(A))}\ar@<-0.5ex>[r]_-{\beta'}&{\Sb}\ar@<-0.5ex>[l]_-{\alpha'}\ar@<-0.5ex>[r]_-{\varphi^{-1}}&{\M(\SB(A))}\ar@<-0.5ex>[l]_-{\varphi}.
}
$$
Consider the cycles $\alpha=\alpha'\circ \varphi$ and $\beta=\varphi^{-1}\circ \beta'$.
Clearly $\alpha\circ\beta=p$.
Consider this equality over the separable closure $k_s$.
We have two rational cycles
$$
\alpha_s=\sum_{i=0}^4 H^i\times a_{4-i} \in\CH^4(\pl^4\times \Gr(2,5)) \text{ and }
$$
$$
\beta_s=\sum_{j=0}^6 b_{6-j}\times H^j \in\CH^6(\Gr(2,5)\times\pl^4)
$$
for some elements $a_i\in \CH^i(\Gr(2,5))$ and $b_j\in \CH^j(\Gr(2,5))$.
Their composition can be written as
$$
\alpha_s\circ\beta_s=\sum_{i=0}^4 b_{6-i}\times a_i=p_s.
$$
Therefore,
$$\alpha_s=\varepsilon_11\times g_4+\varepsilon_2H\times g_3+\varepsilon_3H^2\times (\sigma_2-2g_2)+\varepsilon_4H^3\times\sigma_1+H^4\times 1 \text{ and }
$$
$$
\beta_s=pt\times 1+\varepsilon_4g_5\times H+\varepsilon_3(g_4+3h_4)\times H^2+\varepsilon_2(2\sigma_3+g_3)\times H^3+\varepsilon_1(3\sigma_2+g_2)\times H^4,
$$
where $\varepsilon_i=\pm 1$.
Now let $c$ denotes the rational cycle
$(\rho^2)^t-\alpha_s$ and $b$ the rational cycle
$\rho^3-\beta_s$.
Then the cycle
$$3(b\circ c)=_5 H\times H^3+H^3\times H+3(1-\varepsilon_3)^2H^2\times H^2+3(1-\varepsilon_1)^21\times H^4=d$$
is rational in $\CH^4(\pl^4\times\pl^4)$. Now $d^{\circ 4}$
gives a non-trivial rational projector on $\pl^4$, i.e.,
the motive $\M(\SB(A))$
is decomposable. This leads to a contradiction
with indecomposability of $\M(\SB(A))$ over $\zz$, and finishes the proof of the theorem.

\section{Groups of types $\B_n$ and $\C_n$}\label{Bn}

The goal of the present section is to prove Theorems~\ref{mainB} and \ref{mainC}.

\begin{ntt}
Let $G$ be an adjoint group of type $\B_n$.
From now on we suppose that the characteristic of $k$ is not $2$.
It is known that $G={\mathrm O}^+(V,q)$,
where $(V,q)$ is a regular quadratic space of dimension $2n+1$
and projective $G$-homogeneous varieties can be described as flags
of $q$-totally
isotropic subspaces
$$
X(d_1,\ldots,d_k)=\{V_1\subset\ldots\subset V_k\subset V \mid \dim V_i=d_i\}.
$$
\end{ntt}

\begin{ntt}
The value of the functor corresponding to the variety
$X(d_1,\ldots,d_k)$ at a $k$-algebra $R$
equals the set of all flags
$V_1\subset\ldots\subset V_k$, where
$V_i$ is a $q_R$-totally isotropic direct summand of $V_R$ of rank $d_i$.
\end{ntt}

For convenience we set $d_0=0$, $V_0=0$.

\begin{ntt}
On the scheme $X=X(d_1,\ldots,d_k)$
there are ``tautological'' vector bundles $\V_i$ of ranks $d_i$.
The value of $\V_i$ on an $X$-algebra $(R,x)$ is $V_i$,
where $x=(V_1,\ldots,V_k)$.
For every $m$ there exists an obvious morphism
\begin{align*}
X(d_1,\ldots,d_k)&\to X(d_1,\ldots,\Hat d_m,\ldots,d_k)\\
(V_1,\ldots,V_k)&\mapsto (V_1,\ldots,\Hat V_m,\ldots,V_k)
\end{align*}
which makes $X(d_1,\ldots,d_k)$ into a $X(d_1,\ldots,\hat{d}_m,\ldots,d_k)$-scheme.
\end{ntt}

\begin{lem}\label{GrB}
Denote $X(d_1,\ldots,d_k)$ by $Y$ and
$X(d_1,\ldots,\Hat d_m,\ldots,d_k))$ by $X$.
Suppose that $m<k$. Then $Y$ as a scheme over $X$
can be identified with the Grassmann bundle
$Z=\Gr(d_m-d_{m-1},\V_{m+1}/\V_{m-1})$ over $X$.
\end{lem}

\begin{proof}
Fix an $X$-algebra $(R,x)$,
where $x=(V_1,\ldots,\Hat V_m,\ldots,V_k)$.
We define a natural bijection between the fibers over the point $x$ of $Y$ and $Z$
as follows.
The fiber of $Y$ over $x$
can be identified with the set of all direct summands
$V_m$ of $V_R$ of rank $d_m$ such that $V_{m-1}\subset V_m \subset V_{m+1}$
(note that $V_m$ is automatically $q_R$-isotropic since $V_{m+1}$ is so).
This fiber is clearly isomorphic to the fiber of $Z$ over $x$
which is the set of all direct summands of
$(\V_{m+1}/\V_{m-1})_x=V_{m+1}/V_{m-1}$ of rank $d_m$.
\end{proof}

\begin{ntt}[Proof of Theorem~\ref{mainB}]\label{proofmainB}
Apply Lemma~\ref{GrB} to the varieties $Y=X(d_1,\ldots,d_k)$ and
$X=X(d_1,\ldots,\hat{d}_m,\ldots,d_k)$. We obtain that $Y$ is a Grassmann
bundle over $X$. To finish the proof apply Proposition~\ref{GrBundl}.
\end{ntt}

\begin{ntt}
Let $G$ be an adjoint group of type $\C_n$ over $k$.
It is known that $G=\Aut(A,\sigma)$,
where $A$ is a central simple algebra of degree $2n$
with an involution $\sigma$ of symplectic type on $A$, and
projective $G$-homogeneous varieties can be described as flags
of (right) ideals of $A$
$$
X(d_1,\ldots,d_k)=\{I_1\subset\ldots\subset I_k\subset A\mid I_i\subseteq I_i^\bot,\ \rdim I_i=d_i\}.
$$
Here $I^\bot=\{x\in A\mid\sigma(x)I=0\}$ is a right ideal
of reduced dimension $2n-\rdim I$.
\end{ntt}

\begin{ntt}\label{condC}
The value of the functor corresponding to the variety
$X(d_1,\ldots,d_k)$ at a $k$-algebra $R$ equals
to the set of all flags
$I_1\subset\ldots\subset I_k$ of right ideals of
$A_R=A\otimes_k R$ having the following properties
\begin{itemize}
\item the injection of $A_R$-modules $I_i\hookrightarrow A_R$ splits;
\item $I_i\subseteq I_i^\bot$;
\item $\rdim I_i=d_i$.
\end{itemize}
\end{ntt}

For convenience we set $I_0=0$.

\begin{ntt}
On the scheme $X=X(d_1,\ldots,d_k)$
there are ``tautological'' vector bundles $\J_i$ of ranks
$2nd_i$ and their ``orthogonal complements''
$\J_i^\bot$ of rank $2n(2n-d_i)$.
The value of $\J_i$ (resp. $\J_i^\bot$) on an $X$-algebra $(R,x)$,
where $x=(I_1,\ldots,I_k)$, is $I_i$ (resp. $I_i^\bot$)
considered as a projective $R$-module.
The bundles $\J_i$ and $\J_i^\bot$
also have structures of right $A_X$-modules,
where $A_X$ is a constant sheaf of algebras on $X$ determined by $A$.
There exists an obvious morphism
\begin{align*}
X(d_1,\ldots,d_k)&\to X(d_1,\ldots,d_{k-1})\\
(I_1,\ldots,I_k)&\mapsto (I_1,\ldots,I_{k-1}),
\end{align*}
which makes $X(d_1,\ldots,d_k)$ into a $X(d_1,\ldots,d_{k-1})$-scheme.
\end{ntt}

\begin{lem}\label{grC}
Denote $X(d_1,\ldots,d_k)$ by $Y$ and
$X(d_1,\ldots,d_{k-1})$ by $X$.
Suppose that $d_k=d_{k-1}+1$ and there exists a vector bundle
$\E$ over $X$ such that $A_X\simeq\End_{\OO_X}(\E)$.
Consider the vector bundle
$$
\F=\J_{k-1}^\bot\E/\J_{k-1}\E=\J_{k-1}^\bot/\J_{k-1}\otimes_{A_X}\E
$$
of rank $2(n-d_{k-1})$. Then $Y$ as a scheme over $X$
can be identified with the projective bundle $Z=\pl(\F)=\Gr(1,\F)$ over $X$.
\end{lem}

\begin{proof}
Fix an $X$-algebra $(R,x)$, where $x=(I_1,\ldots,I_{k-1})$.
We define a natural bijection between the fibers over the point $x$ of $Y$ and $Z$.
The fiber of $Y$ can be identified with the set of all
ideals $I_k$ containing $I_{k-1}$ and satisfying the conditions \ref{condC}.
The fiber of $Z$ is the set of all direct summands of $\F_x=\F(R,x)$ of rank $1$.

The involution $\sigma$ induces an isomorphism
$h\colon\E_x\otimes \mathcal{L}\to\E_x^*$ for some invertible $R$-module $\mathcal{L}$
(see \cite[Lemma~III.8.2.2]{Knus}) such that
\begin{align*}
\sigma(f)\otimes 1&=h^{-1}f^*h\text{ for all }f\in A\\
h^*\can\otimes 1&=-h
\end{align*}
where $\can:\E_x\to\E_x^{**}$ is the canonical isomorphism.

Let $U_1$ and $U_2$ be direct summands of $\E_x$.
We write $U_2\subseteq U_1^\bot$
if $h(u\otimes l)(v)=0$ for all $u\in U_1$, $v\in U_2$, $l\in \mathcal{L}$.
We call a direct summand $U$ of $\E_x$ \emph{totally isotropic}
if $U\subseteq U^\bot$.
Note that any direct summand of rank $1$ is totally isotropic (it can be proved easily using localization).

Define $\Phi$ and $\Psi$ as in the proof of Theorem~\ref{thmA}.
Direct computations show that $I_1\subseteq I_2^\bot$ if and only if $\Phi(I_1)\subseteq\Phi(I_2)^\bot$.

So the fiber of $Y$ over $x$ is naturally isomorphic
to the set of all totally isotropic direct summands
$U_k$ of $\E_x$ of rank $d_k$ containing $U_{k-1}=\Phi(I_k)$.
One can represent $U_k$ as the direct sum
$U_{k-1}\oplus U$ where $U$ is a direct summand of rank $1$
(since $d_k=d_{k-1}+1$).
This $U$ is totally isotropic and, therefore,
$U_k$ is totally isotropic if and only if $U_k\subseteq U_{k-1}^\bot$.
Hence the set of all $U_{k-1}$ is naturally isomorphic
to the set of all direct summands of $\Phi(I_{k-1}^\bot)$ of rank $d_k$
containing $\Phi(I_{k-1})$. The latter can be identified with
$\mathbb P(\F_x)$. This finishes the proof of the lemma.
\end{proof}

We are now ready to prove Theorem~\ref{mainC}.

\begin{ntt}[Proof of Theorem~\ref{mainC}]\label{proofmainC}
Consider the flag varieties $Y=X(d_1,\ldots,d_k)$ and
$X(d_1,\ldots,d_{k-1})$. Since $\ind A=2^r$ for some $r$ and there is an odd $d_i$,
we have $\gcd(\ind(A),d_1,\ldots,d_{k-1})=1$.
By Lemma~\ref{propBundle} there exists a bundle $\mathcal{E}$ over $X$ such that
$A_X=\End_{\OO_X}(\mathcal{E})$. Applying Lemma~\ref{grC} to the varieties $X$, $Y$
and the bundle $\mathcal{E}$, we obtain that $Y$ is a projective bundle over $X$.
Now we use Proposition~\ref{GrBundl} to finish the proof.
\end{ntt}

\section{Groups of types $\G$ and $\FF$}\label{except}

This section is devoted to proofs of Theorems~\ref{mainG} and \ref{mainF}.

\begin{ntt}
Let $G$ be a group of type $\G$. We suppose that characteristic of $k$ is not $2$. It is known that $G=\Aut(C)$ where $C$ is a Cayley algebra over $k$. By \emph{$i$-space} where $i=1,2$ we mean an $i$-dimensional subspace $V_i$ of $C$ such that $uv=0$ for every $u,v\in V_i$.

The only flag variety corresponding to a non-maximal parabolic is the full flag variety $X(1,2)$ which is described as follows (see \cite{Bo03}):
$$
X(1,2)=\{V_1\subset V_2\mid V_i\text{ is a $i$-subspace of }C,\ \dim V_i=i\}.
$$

Similarly one can describe the homogeneous flag variety corresponding to the maximal parabolic
$$
X(2)=\{V\mid V\text{ is a $2$-subspace of }C,\ \dim V=2\}.
$$
\end{ntt}

\begin{ntt}
Let $R$ be a $k$-algebra. By an $i$-submodule in $C_R=C\otimes_k R$ we mean a direct summand $V_i$ of $C_R$ of rank $i$ such that $uv=0$ for every two elements $u,v\in V_i$. The value of the functor corresponding to the variety $X(1,2)$ (respectively $X(2)$) at a $k$-algebra $R$ equals the set of all flags $V_1\subset V_2$ (respectively submodules $V_2$) where $V_i$ is an $i$-submodule of $C_R$.
\end{ntt}

\begin{ntt}
On the scheme $X=X(2)$ there is a ``tautological'' vector bundle $\V$ of rank
$2$. The value of $\V$ on an $X$-algebra $(R,x)$ is $V$, where $x=V$.

There exists an obvious morphism
\begin{align*}
X(1,2)&\to X(2)\\
(V_1,V_2)&\mapsto V_2
\end{align*}
which makes $X(1,2)$ into a $X(2)$-scheme.
\end{ntt}

\begin{lem}\label{lemG}
$X(1,2)$ as a scheme over $X(2)$ can be identified with the projective bundle $\mathbb P(\V)=\Gr(1,\V)$ over $X(2)$.
\end{lem}
\begin{proof}
The proof goes as in $\B_n$-case (note that if $V_2$ is a $2$-submodule then each of its direct summands of rank $1$ is a $1$-submodule).
\end{proof}

\begin{ntt}[Proof of Theorem~\ref{mainG}]\label{proofmainG}
Apply Lemma~\ref{lemG} and Proposition~\ref{GrBundl}.
\end{ntt}

\begin{ntt}
Let $G$ be a group of type $\FF$. We suppose that characteristic of $k$ is not $2,3$. It is known that $G=\Aut(J)$ where $J$ is an exceptional Jordan algebra of dimension $27$ over $k$. Set $I=\{1,2,3,6\}$. By \emph{$i$-space} where $i\in I$ we mean an $i$-dimensional subspace $V_i$ of $J$ such that every $u,v\in V_i$ satisfy the following condition:
$$
\tr(u)=0,\ u\times v=0,\text{ and if }i<6\text{ then }u(va)=v(ua)\text{ for all }a\in J.
$$

It is known that projective $G$-homogeneous varieties are parametrized by sequences of numbers $d_1<\ldots <d_k$ from $I$ and can be described as follows:
$$
X(d_1,\ldots,d_k)=\{V_1\subset\ldots\subset V_k\mid V_i\text{ is a $d_i$-subspace of }A,\ \dim V_i=d_i\}.
$$
\end{ntt}

\begin{ntt}
Let $R$ be a $k$-algebra. By an $i$-submodule in $J_R=J\otimes_k R$ we mean a direct summand $V_i$ of $J_R$ of rank $i$ such that every two elements $u,v\in V_i$ satisfy the conditions above. The value of the functor corresponding to the variety $X(d_1,\ldots,d_k)$ at a $k$-algebra $R$ equals the set of all flags
$V_1\subset\ldots\subset V_k$ where $V_i$ is
a $d_i$-submodule of $J_R$.

For convenience we set $d_0=0$, $V_0=0$.
\end{ntt}

\begin{ntt}
On the scheme $X=X(d_1,\ldots,d_k)$ there are ``tautological'' vector bundles $\V_i$ of rank
$d_i$. The value of $\V_i$ on an $X$-algebra $(R,x)$ is $V_i$, where $x=(V_1,\ldots,V_k)$.

There exists an obvious morphism
\begin{align*}
X(d_1,\ldots,d_k)&\to X(d_1,\ldots,\Hat d_m,\ldots,d_k)\\
(V_1,\ldots,V_k)&\mapsto (V_1,\ldots,\Hat V_m,\ldots,V_k)
\end{align*}
which makes $X(d_1,\ldots,d_k)$ into a $X(d_1,\ldots,\Hat d_m,\ldots,d_k)$-scheme.
\end{ntt}

\begin{lem}\label{lemF}
Denote $X(d_1,\ldots,d_k)$ by $Y$ and $X(d_1,\ldots,\Hat d_m,\ldots,d_k)$ by $X$. Suppose that $m<k$ and either $d_{m+1}<6$ or $d_m=1$. Then $Y$ as a scheme over $X$ can be identified with the Grassmann bundle $Z=\Gr(d_m-d_{m-1},\V_{m+1}/\V_{m-1})$ over $X$.
\end{lem}
\begin{proof}
The proof goes as in $\B_n$-case (note that under our restrictions if $V_{m+1}$ is a $d_{m+1}$-submodule then each of its direct summands of rank $d_m$ is a $d_m$-submodule).
\end{proof}

\begin{ntt}[Proof of Theorem~\ref{mainF}]\label{proofmainF}
Apply Lemma~\ref{lemF} to the varieties $Y=X(d_1,\ldots,d_k)$ and
$X=X(d_1,\ldots,\hat{d}_m,\ldots,d_k)$. We obtain that $Y$ is a Grassmann
bundle over $X$. To finish the proof apply Proposition~\ref{GrBundl}.
\end{ntt}

\paragraph{Acknowledgements}
We are sincerely grateful to N.~Karpenko and M.~Rost for the comments
concerning rational cycles and Severi-Brauer varieties.

The first and the last authors would like to thank U.~Rehmann
and RTN-project for hospitality and support.
The last author is also grateful to Alexander von Humboldt foundation
for support.

The second and the third authors highly appreciate A.~Bak's hospitality and sympathy with their activity.
The second author also would like to thank DAAD {\it Jahresstipendium} program and INTAS project
03-51-3251.

\bibliographystyle{chicago}

\end{document}